\documentclass[a4paper,11pt,twoside]{article}
\usepackage{graphicx}
\usepackage{amsmath}
\usepackage{amssymb}
\newtheorem{theorem}{Theorem}[section]
\newtheorem{corollary}[theorem]{Corollary}

\newtheorem{proposition}[theorem]{Proposition}

\newtheorem{definition}[theorem]{Definition}

\def\bkR{{\rm I\kern-.17em R}}
\def\bkC{{\rm \kern.24em
       \vrule width.05em height1.4ex depth-.05ex
       \kern-.26em C}}
\begin{document}

\title{Symbolic approach to the general quadratic polynomial decomposition
}


\author{
{\^A}ngela Macedo \\
Departamento de Matem{\'a}tica and CMAT, Polo CMAT-UTAD\\
Universidade de Tr{\'a}s-os-Montes e Alto Douro (UTAD)\\
Quinta de Prados, 5000-801 Vila Real\\
e-mail: amacedo@utad.pt,  web:  http://www.utad.pt\\
\\
Teresa A. Mesquita\\
Instituto Superior Polit{\'e}cnico de Viana do Castelo and CMUP\\
ESTG, Av. do Atl{\^a}ntico s/n, 4900-348 Viana do Castelo, Portugal\\
e-mail: teresa.mesquita@fc.up.pt, web: http://www.estg.ipvc.pt\\
\\
Z{\'e}lia da Rocha\\
Departamento de Matem{\'a}tica and CMUP\\
Faculdade de Ci{\^e}ncias da Universidade do Porto\\
Rua do Campo Alegre, 687 4169-007 Porto\\
e-mail: mrdioh@fc.up.pt,  \\
web: http://www.fc.up.pt, http://cmup.fc.up.pt\\
}


\date{\today}

\maketitle

\begin{abstract}
In this work we deal with a symbolic approach to the general quadratic polynomial  decomposition.
By means of a symbolic implementation, we investigate some properties of the components sequences like orthogonality and symmetry. We present some explicit results for a collection of well known orthogonal cases.



\end{abstract}


\section{Introduction}
\label{intro}
The quadratic decomposition (QD) of a symmetric orthogonal polynomial sequence into two non-symmetric polynomial sequences goes back to the work of Szeg{\"o}
\cite{GS:1939},
Carlitz \cite{x} and Chihara \cite{Chihara_64,C:78,Chihara_82,ChiharaChihara_87}.
One example is given by the relationship between Hermite and Laguerre polynomials 
\begin{eqnarray}
H_{2n}(x)=\left( -1\right)^n 2^{2n} n! L_{n}^{\left( -1/2\right)}\left( x^2\right)\\
H_{2n+1}(x)=\left( -1\right)^n 2^{2n+1} n! x L_{n}^{\left( 1/2\right)}\left( x^2\right)
\end{eqnarray}
Since then the type of QD has been generalized in order to decompose a non-symmetric sequence \cite{P:90,Quad II} and to admit a full quadratic mapping leading to the so-called general quadratic decomposition (GQD)
\cite{ATese:2004,A:2010}.
In this work we propose a symbolic approach to the GQD in the sequel of what was done in \cite{T:2012} for the general cubic decomposition \cite{Mar_Mes_Rocha_2011}.

In general, when we decompose a polynomial sequence, we may naturally ask if certain properties of the original sequence are also verified by its components and
reciprocally leading to questions of characterization. If we are dealing with an orthogonal sequence the first fundamental property to consider is the regularity
\cite{ATese:2004,A:2010,P:90,Quad II}. We can also study the same problem for specific properties like semi-classical \cite{Bourras_Mar:2012_1,PM-Tounsi:2012},
 Laguerre-Hahn \cite{Bourras_Mar:2010,Bourras_Mar:2012} or Appell characters \cite{AFL:2008}, among others.
 Also, the decomposition can reveal interesting connections between a sequence and its components and allows a deep study of all families involved as it was the case in \cite{AFL:2011,TA:2015}.
Symmetrization and the application of other polynomial mappings have always been a constructive method for obtaining new orthogonal families \cite{JesusPetro_2010,Mar_San:1991,Mar_San:1992,Mar_Pet:1997}.
Polynomial mappings lead to special kinds of orthogonality measures \cite{Ger_VanAssche_1988} and are related to the so-called
{\em sieved polynomials} \cite{Al_All_As_1984,Ismail_86_3}.

The study of different types of polynomial sequences demands frequently a great amount of algebraic recursive computations which motivates the introduction of symbolic software as a fundamental research tool in the framework of orthogonal polynomials.
The recurrence relations established in this field fit naturally within an algorithmic
approach and the symbolic computation is currently growing towards a consistent partner in theoretical and applied research in this area.
This has been the case in some contributions, we indicate as examples the references
\cite{PZ:08,PZ:13,ZR:2016,ZR:2016_1}, where the computational environment was determinant to give explicitly several theoretical properties regarding, in particular,
semi-classical orthogonal polynomials.
In this work we adopt the same point of view, we exploit the capabilities of symbolic computations in order to study some properties and to explicit closed formulas related to the GQD.

This article is organized as follows. In Section 2 we recall some notation and basic concepts about polynomial and dual sequences, symmetry and orthogonality that will be needed in the sequel. Section 3 is devoted to the GQD, its definition and properties when the original sequence is orthogonal; in particular we have collected from
\cite{A:2010,ATese:2004} the theoretical results and formulas that we have implemented in this work.
In Section 4 we describe step by step the symbolic implementation done in order to compute and investigate some properties of the component sequences.
In Section 5 we begin by considering two study cases of orthogonal sequences, one symmetric and another non-symmetric. We decompose these two examples and give in detail new results about their GQD provided by the symbolic implementation. Next, we study a collection of well known orthogonal sequences consisting in all classical families, the generalized Hermite and the Charlier sequences.
By means of the symbolic implementation, we explicit new closed formulas for the extended coefficients for arbitrary parameters.  If one, two or three of the parameters presented in the GQD vanish, we furnish also the indication of orthogonality
or non-orthogonality and the existing vanishing component sequences.
We finish this article with some conclusions about the implementation and results obtained.

\section{Notation and basic concepts}
Let $\mathcal{P}$ be the vector space of polynomials with coefficients in $\bkC$ and let
$\mathcal{P}^{\prime}$ be its topological dual space. We denote by $\langle w ,p\rangle$ the effect of $w \in\mathcal{P}^{\prime }$ on $p\in\mathcal{P}$. In particular,
$\langle w, x^{n}\rangle:=\nolinebreak\left( w \right) _{n},n\geq 0$, represent the {\it moments} of $w$.
In the following, we will call polynomial sequence (PS) to any sequence
${\{W_{n}\}}_{n \geq 0}$ such that $\deg W_{n}= n,\; \forall n \geq 0$.
We will also call monic polynomial sequence (MPS) to a PS so that all polynomials have leading coefficient equal to one.
If ${\{W_{n}\}}_{n \geq 0}$ is a MPS, there exists a unique sequence
$\{w_n\}_{n\geq 0}$, $w_n\in\mathcal{P}^{\prime}$, called the {\it dual sequence} of $\{W_{n}\}_{n\geq 0}$, such that,
\begin{equation}\label{SucDual}
<w_{n},W_{m}>=\delta_{n,m}\ , \ n,m\ge 0.
\end{equation}
The form $w_0$ is called the {\it canonical form} of  $\{W_{n}\}_{n\geq 0}$; it is normalized, i.e.,
$(w_0)_0=\langle w_{0}, 1 \rangle=1$.
Let us write $W_m(x)$ in the canonical basis
\begin{equation}\label{PnCan}
W_0(x)=1\quad,\quad W_m(x)=\sum_{\nu=0}^{m-1}a_{m,\nu}x^{\nu}+x^m\ ,\ m\geq 1\ .
\end{equation}
Replacing (\ref{PnCan}) in  (\ref{SucDual}), we obtain the moments of the form $w_n$ in terms of the canonical coefficients of $W_m(x)$ \cite[p.144]{ZR:1999}
\begin{equation}\label{MomSucDual}
(w_n)_m=\left\{
\begin{array}{l l}
0\ , & m=0,\ldots,n-1,\\
1\ , & m=n,\\
-\displaystyle\sum_{\nu=n}^{m-1}a_{m,\nu}(w_n)_\nu\ , & m=n+1,n+2,\ldots
\end{array}
\right..
\end{equation}
In particular, taking $n=0$, we get the moments of the canonical form $w_0$
\begin{equation}\label{MomSucCan}
(w_0)_0=1\quad, \quad (w_0)_m=-\sum_{\nu=0}^{m-1}a_{m,\nu}(w_0)_\nu\ ,\ m\geq 1.
\end{equation}

%
\noindent
Furthermore, given a MPS ${\{W_{n}\}}_{n \geq 0}$, applying the euclidian division of $W_{n+2}$ by $W_{n+1}$, we can always define complex sequences, ${\{\beta_{n}\}}_{n \geq 0}$ and
$\{\chi_{n,\nu}\}_{0 \leq \nu \leq n,\; n \geq 0},$ such that
\begin{eqnarray}
&&W_{0}(x)=1, \;\; W_{1}(x)=x-\beta_{0},\label{divisao_ci}\\
 &&W_{n+2}(x)=(x-\beta_{n+1})W_{n+1}(x)-\sum_{\nu=0}^{n}\chi_{n,\nu}W_{\nu}(x).\label{divisao}
\end{eqnarray}
This relation is usually called the structure relation of  ${\{W_{n}\}}_{n \geq 0}$, and ${\{\beta_{n}\}}_{n \geq 0}$ and
$\{\chi_{n,\nu}\}_{0 \leq \nu \leq n,\; n \geq 0}$ are called the structure coefficients. Their relationship with the dual sequence is given by \cite{P:94,P:91}
\begin{equation}\label{SC_dualSeq}
\beta_n=<w_n,xW_n>\ ,\  n\geq 0\ \ ;\ \
\chi_{n,\nu}=<w_\nu,xW_{n+1}>\ , \ 0\leq \nu\leq n\ ,\ n\geq 0\ .
\end{equation}
\medskip

A PS $\{W_{n}(x)\}_{n \geq 0}$ is said to be {\it symmetric} if and only if $W_{n}(-x)=(-1)^{n}W_{n}(x),$ $\; n \geq 0$.
Moreover, a MPS ${\{W_{n}\}}_{n\geq 0}$ is symmetric if and only if $(w_{2n})_{2\nu+1}=0$ and $(w_{2n+1})_{2\nu}=0$, $\nu\geq n$, $n\geq 0$ if and only if \cite{P:91}
\begin{equation}\label{}
\beta_{n}=0\ ;\ \chi_{2n+1,2\nu}=0, 0 \leq \nu \leq n, n \geq 0\ ;\
\chi_{2n,2\nu+1}=0,0 \leq \nu \leq n-1,n \geq 1.
\end{equation}

\begin{definition} \cite{P:91,P:94} \label{orthogonality definition}
A PS ${\{W_{n}\}}_{n \geq 0}$ is (regularly) orthogonal with respect to the form $w$ if and only if it fulfils
\begin{eqnarray}
&&\label{ortogonal} <w,W_{n}W_{m}>=0,\:\: n \neq m,\:\:\:\:\: n, m \geq 0,\\
&&\label{ortogonal regular} <w,W_{n}^{2}> \neq 0, \: n \geq 0.
\end{eqnarray}
It is easy to prove that necessarily $w = (w)_0 \,w_{0}$ \cite{P:91}. ${\{W_{n}\}}_{n \geq 0}$ and $w$ can always be taken normalized in such a way ${\{W_{n}\}}_{n \geq 0}$ is MOPS and $w=w_0$.
\end{definition}

Let ${\{W_{n}\}}_{n \geq 0}$ be an orthogonal MPS (shortly indicated as MOPS). It is well known that the structure relation (\ref{divisao_ci})-(\ref{divisao})
becomes the following second order recurrence relation:
\begin{eqnarray}
 \label{RSecondOrder-I} &&  W_{0}(x)=1,\;\;W_{1}(x)=x-\beta_{0},\\
\label{RSecondOrder} && W_{n+2}(x)=(x-\beta_{n+1})W_{n+1}(x)-\gamma_{n+1}W_{n}(x),\;\;\; n \geq 0,
\end{eqnarray}
where $\gamma_{n+1}=\chi_{n,n}\neq 0,\;\; n \geq 0$ and
\begin{equation}
\beta_n=\frac{<w_0,xP_n^2>}{<w_0,P_n^2>}\ ,\  n\geq 0\ ;\
\gamma_{n+1}=\frac{<w_0,P_{n+1}^2>}{<w_0,P_n^2>}\ , \ 0\leq \nu\leq n\ ,\ n\geq 0\ .
\end{equation}
Note that the regularity conditions (\ref{ortogonal regular}) are fulfilled if and only if $\gamma_{n+1} \neq 0,\:\: n \geq 0$. The coefficients $\beta_{n}$ and $\gamma_{n+1} $ are called the recurrence coefficients of ${\{W_{n}\}}_{n \geq 0}$.

A MOPS $\{W_n\}_{n\geq 0}$ with respect to $w_0$ is real if and only if $\beta_{n}\in \bkR$ and $\gamma_{n+1}\in \bkR\setminus\{0\}$. These conditions are equivalent to $(w_0)_n\in \bkR$,
$n\geq 0$, i.e., $w_0$ is real. Furthermore, $\beta_{n}\in \bkR$ and $\gamma_{n+1}>0$, $n\geq 0$ if and only if $w_0$ is {\em positive definite}
if and only if  $(w_0)_n\in \bkR$,
$n\geq 0$, and $(w)_{2n}>0$, $n\geq 0$ \cite{x,P:91}.  A MOPS $\{W_n\}_{n\geq 0}$ is symmetric if and only if $\beta_n=0$, $n\geq 0$ if and only if $(w_0)_{2n+1}=0$, $n\geq 0$, i.e., $w_0$ is symmetric \cite{x,P:91}.

Let ${\{W_{n}\}}_{n \geq 0}$ be a MPS.  Given $A \in \bkC \backslash \{0\}$ and $B \in \bkC $, let us define the MPS provided by the 
shifting mapping $x \mapsto Ax+B$, as follows
\begin{equation}\label{shifted sequence}
\tilde{W}_{n}(x)=A^{-n}W_{n}(Ax+B)\ , \quad n \geq 0\ .
\end{equation}
If ${\{W_{n}\}}_{n \geq 0}$ is a MOPS, then the MPS defined by (\ref{shifted sequence})
is also a MOPS
and its recurrence coefficients are \cite{C:78,variations}
\begin{equation}\label{shifted-rec-coeff}
\widetilde{\beta}_{n}=\frac{\beta_{n}-B}{A}\quad ,\quad \widetilde{\gamma}_{n+1}=\frac{\gamma_{n+1}}{A^{2}}\ ,\quad n \geq 0.\end{equation}

Considering a MOPS and thus fulfilling the recurrence (\ref{RSecondOrder-I})-(\ref{RSecondOrder}),
we may define the associated sequence of order $r$ with $r\geq 0$, a co-recursive sequence,
and a wider perturbation of order $r$, with $r\geq 1$ (which includes a strictly co-recursive situation) as the next definitions explain.

\begin{definition}\cite{P:91}\label{def-ass-r}
Let $\{W_{n}\}_{n\ge 0}$ be a MOPS satisfying (\ref{RSecondOrder-I})-(\ref{RSecondOrder}). The  MOPS $\{W_{n}^{(r)}\}_{n\ge 0}$ defined by
\begin{equation}\label{(1.4.3.3)}
\left\{\begin{array}{l}
W_{0}^{(r)}(x)=1\quad,\quad W_{1}^{(r)}(x)=x-\beta_{0}^{(r)}\;,\\
W_{n+2}^{(r)}(x)=(x-\beta_{n+1}^{(r)})W_{n+1}^{(r)}(x)-\gamma_{n+1}^{(r)}W_{n}^{(r)}(x)\ ,\quad n\;,\;r\geq 0\;,
\end{array}
\right.
\end{equation}
where
\begin{equation}\label{(1.4.3.4)}
\beta_{n}^{(r)}=\beta_{n+r}\quad,\quad\gamma_{n+1}^{(r)}=\gamma_{n+1+r}\ ,\quad n\;,\;r\ge 0\;.
\end{equation}
is called the associated sequence of order $r$ of $\{W_{n}\}_{n\ge 0}$.
\end{definition}
\noindent

\begin{definition}\cite{P:91}\label{defn-co-rec}
Let ${\{W_{n}\}}_{n \geq 0}$ be a MOPS satisfying (\ref{RSecondOrder-I})-(\ref{RSecondOrder}).
Given $\mu \in \bkC$, the MOPS ${\{W_{n}(\mu; x)\}}_{n \geq 0}$ defined by
\begin{equation}\left\{ \begin{array}{l}
W_{0}(\mu;x)=1\quad,\quad W_{1}(\mu; x)=x-\beta_{0}-\mu\;,\\
W_{n+2}(\mu; x)=(x-\beta_{n+1})W_{n+1}(\mu; x)-\gamma_{n+1}W_{n}(\mu; x)\ ,\quad n\ge 0\;,
\end{array}
\right.
\end{equation}
is called a co-recursive sequence of ${\{W_{n}\}}_{n \geq 0}$.
\end{definition}

\begin{definition}\cite{P:91}\label{defn r-perturbed}
Let ${\{W_{n}\}}_{n \geq 0}$ be a MOPS satisfying (\ref{RSecondOrder-I})-(\ref{RSecondOrder}) and let $r \geq 1$.
\newline
Given $\mu_{0} \in \bkC$,  $\mu=(\mu_{1},\ldots,  \mu_{r})\: \in \bkC^{r}$ and $\lambda=(\lambda_{1},\ldots,  \lambda_{r})\:
\in \Big(\bkC\setminus \{0\}\Big)^{r}$, where either $ \mu_{r} \neq 0$ or $ \lambda_{r} \neq 1$,
\begin{equation}\label{}
\left\{\begin{array}{l}
\tilde{W}_{0}(x)=1,\quad\tilde{W}_{1}(x)=x-\tilde{\beta}_{0}, \\
\tilde{W}_{n+2}(x)=(x-\tilde{\beta}_{n+1})\tilde{W}_{n+1}(x)-\tilde{\gamma}_{n+1}\tilde{W}_{n}(x)\ ,\quad n \geq 0\ ,
\end{array}
\right.
\end{equation}
with
$$\begin{array}{ll}
\tilde{\beta}_{0}=\beta_{0}+\mu_{0} &,\quad \\
\tilde{\beta}_{n}=\beta_{n}+\mu_{n} &,\quad \tilde{\gamma}_{n}=\lambda_{n}\gamma_{n} \ ,\quad 1\leq n \leq r\ ,\\
\tilde{\beta}_{n}=\beta_{n} & , \quad  \tilde{\gamma}_{n}=\gamma_{n} \ \ \ \ ,\quad n \geq r+1\ ,
\end{array}$$
is called a perturbed sequence of order $r$ of ${\{W_{n}\}}_{n \geq 0}$ and is denoted by
$${\Big\{W_{n}\Big(\mu_{0};  \begin{array}{c}
\mu\\
\lambda\\
\end{array}; r ; x
 \Big)\Big\}}_{n \geq 0}.$$
\end{definition}

\section{General Quadratic Decomposition}
In \cite{ATese:2004,A:2010} it was introduced the general quadratic
decomposition as follows.
Given a MPS $\{W_n\}_{n\geq 0}$, and fixing a monic quadratic polynomial mapping
\begin{equation}
\:\omega(x)=x^2+px+q\;,\:\:p,q\in \bkC
\end{equation}
and a constant $a\in \bkC,$ it is always possible to associate with it, two MPSs  $\:\left\{P_n\right\}_{n\geq 0}\:$ and $\:\left\{R_n\right\}_{n\geq 0}\:$
and two other sequences of polynomials  $\:\left\{a_n\right\}_{n\geq 0}\:$
and $\:\left\{b_n\right\}_{n\geq 0}\:$ through the unique decomposition
\begin{eqnarray}
W_{2n}(x)&=&P_n(\omega (x))+(x-a) a_{n-1}(\omega(x)) \label{eq1}\\
W_{2n+1}(x)&=&b_n(\omega (x))+(x-a) R_{n}(\omega (x)),\;\; n \geq 0, \label{eq2}
\end{eqnarray}
where  $\:\deg
\left( a_n(x)\right) \leq n\:,\:\:\deg \left( b_n(x) \right) \leq n\:,$ and $a_{-1}(x)=0\:.\:$
$\:\left\{P_n\right\}_{n\geq 0}\:$ and $\:\left\{R_n\right\}_{n\geq 0}\:$ are called the principal component sequences and $\:\left\{a_n\right\}_{n\geq 0}\:$
and $\:\left\{b_n\right\}_{n\geq 0}\:$ are said the secondary ones. Let us note the dual sequences corresponding to $\:\left\{P_n\right\}_{n\geq 0}\:$ and $\:\left\{R_n\right\}_{n\geq 0}\:$ by $(u_n)_{n\geq0}$ and $(v_n)_{n\geq0}$, respectively.

As the principal components are MPS, then given $\{W_n\}_{n\geq 0}$, there exists two unique tables of coefficients
\cite{ATese:2004,A:2010,P:91}, $\:\left(\lambda_\nu^n\right)\:,\:
\left(\theta_\nu^n\right)\:,\: 0\leq \nu\leq n\:,\: n\geq 0\:$,
such that:
\begin{equation}
a_{n}(x)=\sum_{\nu =0}^n \lambda_\nu^n R_\nu (x)\:\:\mbox{and}\:\:
b_{n}(x)=\sum_{\nu =0}^n \theta_\nu^n P_\nu (x)\:,\:n\geq 0\:.
\label{eq3}
\end{equation}
Reciprocally, given $\{P_n\}_{n\geq 0}$ and $\{R_n\}_{n\geq 0}$ and the coefficients $\:\left(\lambda_\nu^n\right)\:,\:
\left(\theta_\nu^n\right)\:$, $0\leq \nu\leq n\:,\: n\geq 0\:$, the sequence
$\{W_n\}_{n\geq 0}$ is determined uniquely by (\ref{eq1}) and (\ref{eq2}).

An important particular case occurs when all parameters $p$, $q$ and $a$ are equal to zero. In this situation, if $\{W_n\}_{n\geq 0}$ is symmetric,
then the secondary components $\:\left\{a_n\right\}_{n\geq 0}\:$ and $\:\left\{b_n\right\}_{n\geq 0}\:$ vanish. This is the case in the  example mentioned in the introduction gathering the Hermite and the Laguerre sequences.
When the secondary components vanish, $\{W_n\}_{n\geq 0}$
is called $\omega$-symmetric \cite{ATese:2004,A:2010}. We notice that the QD of a non-symmetric sequence with vanishing parameters is extensively studied in \cite{P:90,Quad II}.

We will now state the theoretical results needed in the symbolic implementation presented in Section \ref{sectionSIGQD} closely following \cite{ATese:2004,A:2010}.

\begin{proposition}\cite{ATese:2004,A:2010}
\label{prop_nova}
A MPS $\left\{W_n\right\}_{n\geq 0}\:$
defined by (\ref{divisao_ci})-(\ref{divisao}) fulfils (\ref{eq1})-(\ref{eq2}) if and only if the following recurrence relations hold.
\begin{eqnarray}
\label{pascaleq5.0} && \; b_0(x)=a-\beta_0 \notag \\
\label{pascaleq5} &&\; P_{n+1}(x)=-\displaystyle\sum_{\nu=0}^n\chi_{2n,2\nu}P_\nu(x)+\left(x-\omega(a)\right)R_n(x)
+\left(a-\beta_{2n+1}\right)b_n(x)\notag \\
\nonumber &&\quad \quad \quad \quad -\displaystyle\sum_{\nu=0}^{n-1}\chi_{2n,2\nu+1}b_\nu(x) \\
\label{pascaleq6} &&\; a_n(x)=-\displaystyle\sum_{\nu=0}^n\chi_{2n,2\nu}a_{\nu-1}(x)-\left(a+p+\beta_{2n+1}\right)R_n(x)+b_n(x)  \notag\\
\nonumber  &&\quad \quad \quad \quad -\displaystyle\sum_{\nu=0}^{n-1}\chi_{2n,2\nu+1}R_\nu(x)\\
\label{pascaleq7} && \; b_{n+1}(x)=-\displaystyle\sum_{\nu=0}^n\chi_{2n+1,2\nu+1}b_\nu(x)+\left(a-\beta_{2n+2}\right)P_{n+1}(x) \notag\\
\nonumber  &&\quad \quad \quad \quad +
\left(x-\omega(a)\right)a_n(x) -\displaystyle\sum_{\nu=0}^n\chi_{2n+1,2\nu}P_\nu(x)\notag\\
\label{pascaleq8} && \; R_{n+1}(x)=-\displaystyle\sum_{\nu=0}^n\chi_{2n+1,2\nu+1}R_\nu(x)+P_{n+1}(x)-\left(a+p+\beta_{2n+2}\right)a_n(x)\notag \\
\nonumber  &&\quad \quad \quad \quad  -\displaystyle\sum_{\nu=0}^n\chi_{2n+1,2\nu}a_{\nu-1}(x)
\end{eqnarray}
\end{proposition}

We point out that among the component sequences the principal ones $\:\left\{P_n\right\}_{n\geq 0}\:$ and $\:\left\{R_n\right\}_{n\geq 0}\,\:$ are MPSs,
and thus we may ask if they are  orthogonal when the given sequence $\:\left\{W_n\right\}_{n\geq 0}\:$ is orthogonal.
The answer is negative as we can realize from the next proposition or by some computational experiments.

\begin{proposition}\cite{ATese:2004,A:2010} \label{QD-MOPS}
A MPS $\:\left\{W_n\right\}_{n\geq 0}\:$ defined
by (\ref{eq1})-(\ref{eq2}) is  orthogonal if and only if the following recurrence relations are fulfilled for $n\geq 0$,
\begin{equation}
\begin{array}{l}
a_0(x)=-\left(p+\beta_0+\beta_1\right)\\
b_0(x)=a-\beta_0\\
b_{n+1}(x)=-\gamma_{2n+2}b_{n}(x)+\left(a-\beta_{2n+2}\right)P_{n+1}(x)+\left(x-\omega(a)\right)a_n(x)\\
a_{n+1}(x)=-\gamma_{2n+3}a_{n}(x)-\left(a+p+\beta_{2n+3}\right)R_{n+1}(x)+b_{n+1}(x)\\
\end{array}
\label{A_nB_MOPS}
\end{equation}
\begin{equation}\label{RR_ERC_Ps}
\left\{
\begin{array}{l}
P_0(x)=1\:\:,\quad P_1(x)=x-\beta_0^P\ ,\\[1em]
P_{n+2}(x) =\Big(x-\beta_{n+1}^P\Big)P_{n+1}(x)
-\gamma_{n+1}^PP_n(x)
-\varrho_{n+1}^Pb_{n+1}(x)-\rho_{n+1}^Pb_{n}(x)\ , \notag
 \end{array}
 \right.
\end{equation}
with the following extended coefficients of $\{P_n\}_{n\geq 0}$ and
\begin{equation}
\begin{array}{l}
\beta_0^P=\gamma_1+\omega(a)-\left(a-\beta_0\right)\left(a-\beta_1\right) \\ 
\beta_{n+1}^P=\omega(a)+\gamma_{2n+2}+\gamma_{2n+3}-\left(\beta_{2n+2}+a+p\right)\left(a-\beta_{2n+2}\right)\\
\gamma_{n+1}^P=\gamma_{2n+1}\gamma_{2n+2}\\
\varrho_{n+1}^P=\left(\beta_{2n+2}+\beta_{2n+3}+p\right)\\
\rho_{n+1}^P=\gamma_{2n+2}\left(\beta_{2n+1}+\beta_{2n+2}+p\right)
\end{array}
\label{ERC_P}
\end{equation}
\begin{equation}\label{RR_ERC_Rs}
\left\{
\begin{array}{l}
R_0(x)=1\:\:,\quad R_1(x)=x-\beta_0^R\ , \\[1em]
R_{n+2}(x)=\Big(x-\beta_{n+1}^R\Big)R_{n+1}(x)-\gamma_{n+1}^RR_n(x)-\varrho_{n+1}^Ra_{n+1}(x)-\rho_{n+1}^Ra_{n}(x)\ , \notag
  \end{array}
 \right.
\end{equation}
with the following extended coefficients of $\{R_n\}_{n\geq 0}$
\begin{equation}\label{ERC_Rs}
\begin{array}{l}
\beta_0^R=\omega(a)+\gamma_1+\gamma_2-\left(a-\beta_0\right)\left(a-\beta_1\right)-
\left(p+\beta_0+\beta_1\right)\left(\beta_2+a+p\right)\\ 
\beta_{n+1}^R=\omega(a)+\gamma_{2n+3}+\gamma_{2n+4}-\left(\beta_{2n+3}+a+p\right)
\left(a-\beta_{2n+3}\right)\\
\gamma_{n+1}^R=\gamma_{2n+2}\gamma_{2n+3}\\
\varrho_{n+1}^R=\left(\beta_{2n+3}+\beta_{2n+4}+p\right)\\
\rho_{n+1}^R=\gamma_{2n+3}\left(\beta_{2n+2}+\beta_{2n+3}+p\right)\\
\end{array}
\end{equation}
\end{proposition}
This proposition is equivalent to the following one.

\begin{proposition}\cite{ATese:2004}\label{QD-MOPS-anbn}
A MPS $\:\left\{W_n\right\}_{n\geq 0}\:$ defined
by (\ref{eq1})-(\ref{eq2}) is  orthogonal if and only if the following recurrence relations are fulfilled for $n\geq 0$,
\begin{equation}
P_{n+1}(x)=-\gamma_{2n+1}P_n(x)+\left(x-\varpi(a)\right)R_n(x)
+\left(a-\beta_{2n+1}\right)b_n(x)
\label{sist6.c}
\end{equation}
\begin{equation}
R_{n+1}(x)=-\gamma_{2n+2}R_n(x)+P_{n+1}(x)-\left(a+p+\beta_{2n+2}\right)a_n(x)
\label{sist6.d}
\end{equation}
\begin{equation}
\left(a+p+\beta_{2n+1}\right)\left(a-\beta_{2n+2}\right)\neq 0
\label{sist6.e}
\end{equation}
\begin{equation}\left\{
\begin{array}{l}
a_0(x)=-(p+\beta_0+\beta_1)\,,\\[.5em]
\begin{split}
a_{n+1}(x)&=\Big\{(a-\beta_{2n+2})^{-1}(a+p+\beta_{2n+3})(x-\omega(\beta_{2n+2}))\\
&-\gamma_{2n+2}(a+p+\beta_{2n+1})^{-1}(a+p+\beta_{2n+3})
-\gamma_{2n+3}
\Big\}a_{n}(x)\\[.7em]
&
-\gamma_{2n+1}\gamma_{2n+2}(a+p+\beta_{2n+1})^{-1}(a+p+\beta_{2n+3})a_{n-1}(x)\\[.7em]
&
-(a-\beta_{2n+2})^{-1}(p+\beta_{2n+2}+\beta_{2n+3})b_{n+1}(x)\\[.7em]
& -\gamma_{2n+2}(a+p+\beta_{2n+1})^{-1}(a-\beta_{2n+2})^{-1}(a+p+\beta_{2n+3})\times\\
& (p+\beta_{2n+1}+\beta_{2n+2})b_{n}(x)\,,
\end{split}
 \end{array}\right.
 \label{sist6.1}
\end{equation}
\begin{equation}\left\{
\begin{array}{l}
b_0(x)=a-\beta_0\,,\quad
\\
\begin{split}b_1(x)=&\Big\{a-(p+\beta_0+\beta_1+\beta_2)\Big\}(x-\omega(a))-\gamma_1(a-\beta_2)
\\&
-\gamma_2(a-\beta_0)+(a-\beta_0)(a-\beta_1)(a-\beta_2)\,,\\[.5em]
\end{split}\\
[1em]
\begin{split}
b_{n+2}(x)=&\Big\{(a+p+\beta_{2n+3})^{-1}(x-\omega(\beta_{2n+3}))(a-\beta_{2n+4})\\
&
-\gamma_{2n+3}(a-\beta_{2n+2})^{-1}(a-\beta_{2n+4})
-\gamma_{2n+4}
\Big\}b_{n+1}(x)\\[.7em]
&
-\gamma_{2n+2}\gamma_{2n+3}(a-\beta_{2n+2})^{-1}(a-\beta_{2n+4})b_{n}(x)\\[.7em]
&
+(a+p+\beta_{2n+3})^{-1}(p+\beta_{2n+3}+\beta_{2n+4})\times\\
&(x-\omega(a))a_{n+1}(x)\\[.7em]
& +\gamma_{2n+3}(a-\beta_{2n+2})^{-1}(p+\beta_{2n+2}+\beta_{2n+3})\\
&(a+p+\beta_{2n+3})^{-1}(a-\beta_{2n+4})(x-\omega(a))a_{n}(x)\,.
\end{split}
  \end{array}\right.
\label{sist6.2}
\end{equation}
\end{proposition}

\noindent
Relations (\ref{sist6.1}) and (\ref{sist6.2}) involve only the secondary sequences $\{a_n(x)\}_{n\geq0}$ and
$\{b_n(x)\}_{n\geq0}$. They will be crucial in Section \ref{results}, because they allow to do an induction reasoning to prove,
for all $n\geq 0$, some models for secondary sequences inferred from the first polynomials computed by the symbolic implementation of Section \ref{sectionSIGQD}.

\begin{corollary}\label{corolarioOC}
Let $\:\left\{W_n\right\}_{n\geq 0}\:$ be a MOPS. If the extended coefficients $\:\left\{\varrho_n^P\right\}_{n\geq 1}\,,\:$
$\:\left\{\rho_n^P\right\}_{n\geq 0}\:$,
$\:\left\{\varrho_n^R\right\}_{n\geq 1}\:$  and
$\:\left\{\rho_n^R\right\}_{n\geq 0}\:$ are equal to zero, then the MPS
$\:\left\{P_n\right\}_{n\geq 0}\:$ and
$\:\left\{R_n\right\}_{n\geq 0}\:$ are orthogonal.
\end{corollary}
This corollary is obtained directly from (\ref{ERC_P}) and (\ref{ERC_Rs}), because under its hypothesis the principal components will verify recurrence relations of order two of type
(\ref{RSecondOrder-I})-(\ref{RSecondOrder}) and consequently they are orthogonal.
We notice that the reciprocal of this corollary is not true. In the case $\:\left\{W_n\right\}_{n\geq 0}\:$ is a MOPS, the characterization of the orthogonality of the principal components is provided by the next proposition.
\begin{proposition}\cite{ATese:2004,A:2010} \label{prop4}
Let $\:\left\{W_n\right\}_{n\geq 0}\:$ be a MOPS. The MPS  $\:\left\{P_n\right\}_{n\geq 0}\:$
is orthogonal if and only if the coefficients $\:\lambda_\nu^n\:$
and $\:\theta_\nu^n\,,\:$ given in (\ref{eq3}), satisfy:
\begin{eqnarray}
&& \left(\beta_{2n+2}+\beta_{2n+3}+p\right)\theta_\nu^{n+1}+ \left(\beta_{2n+1}+\beta_{2n+2}+p\right)\theta_\nu^{n}=0\ , \notag\\ 
&& \hspace{6.7cm}  0\leq \nu\leq n-1\ ,\  n\geq 1\ , \notag\\
&&\gamma_{n+1}^P:=\gamma_{2n+1}\gamma_{2n+2}+\left(\beta_{2n+2}+\beta_{2n+3}+p\right)\theta_n^{n+1}+\notag\\
&& \hspace{3.75cm}  \gamma_{2n+2}
\left(\beta_{2n+1}+\beta_{2n+2}+p\right)\theta_n^{n}\neq
0\ ,\ n\geq 0\ . \notag 
\end{eqnarray}
Likewise,
 the MPS $\:\left\{R_n\right\}_{n\geq 0}\:$ is
orthogonal if and only if the coefficients $\:\lambda_\nu^n\:$ and
$\:\theta_\nu^n\:$, given in (\ref{eq3}), satisfy:
\begin{eqnarray}
&& \left(\beta_{2n+3}+\beta_{2n+4}+p\right)\lambda_\nu^{n+1}+\gamma_{2n+3}
\left(\beta_{2n+2}+\beta_{2n+3}+p\right)\lambda_\nu^{n}=0\ , \notag\\ 
&& \hspace{6.75cm}  0\leq \nu\leq n-1\ ,\ n\geq 1 \ , \notag\\
&& \gamma_{n+1}^R:=\gamma_{2n+2}\gamma_{2n+3}+\left(\beta_{2n+3}+\beta_{2n+4}+p\right)\lambda_n^{n+1}+ \notag\\
&& \hspace{3.75cm}  \gamma_{2n+3}
\left(\beta_{2n+2}+\beta_{2n+3}+p\right)\lambda_n^{n}\neq 0\ ,\  n\geq 0\ . \notag 
\end{eqnarray}
\end{proposition}

\section{Symbolic implementation of the general quadratic decomposition}\label{sectionSIGQD}

We would like to decompose a MOPS $\{W_n\}_{n\geq 0}$ and to obtain the four component sequences of its GQD in order to study their properties and reveal relationships between them and with $\{W_n\}_{n\geq 0}$.
The original sequence $\{W_n\}_{n\geq 0}$ is given by its recurrence coefficients $\{\beta_{n}\}_{n \geq 0}$ and $\{\gamma_{n+1}\} _{n \geq 0}$. In the majority of examples, the closed formulas of these
coefficients are known for all $n$. If not, the knowledge of their first elements is sufficient for computing the first polynomials of the components.
We fix the parameters of the decomposition, that is, the mapping stated by the polynomial $\omega(x)=x^2+px+q$, given by its coefficients $p$ and $q$, and the zero $a$ of the auxiliary polynomial $(x-a)$ appearing in (\ref{eq1}) and (\ref{eq2}). Also, we fix a positive integer $nmax$, for the maximal degree of the first polynomials we are going to compute.

Established the starting data, we begin by computing the closed formulas, valid for all $n$, to the extended coefficients of the principal components from identities (\ref{ERC_P}) and (\ref{ERC_Rs}) of Proposition \ref{QD-MOPS}. If the last two extended coefficients of $\{P_n\}_{n\geq 0}$ or $\{R_n\}_{n\geq 0}$, i.e., $\varrho_{n+1}^P$, $\rho_{n+1}^P$, $n\geq 0$, or $\varrho_{n+1}^R$, $\rho_{n+1}^R$, $n\geq 0$, vanish, we conclude that the corresponding principal component is orthogonal.
We have also the possibility to compare the recurrence coefficients of the principal components with those of $\{W_n\}_{n\geq 0}$ with the aim to identify any existing relationship between all  these sequences with respect to shifting, association and perturbation transformations.
The first moments of the canonical forms of $\{P_n\}_{n\geq 0}$ and $\{R_n\}_{n\geq 0}$ can always be computed recursively from the canonical coefficients of the
first polynomials using the general procedure stated by (\ref{MomSucCan}). If the last two extended coefficients of $\{P_n\}_{n\geq 0}$ or $\{R_n\}_{n\geq 0}$
are non-zero, the component can or can not be orthogonal. We may inquire about this property by computing the structure coefficients of the first polynomials of the sequence making the euclidian division, or using the direct identities (\ref{SC_dualSeq}) which demand the first moments of the dual sequences. Then, we can analyse the structure coefficients obtained and conclude about the orthogonality of the first polynomials, that occurs if $\chi_{n,\nu}=0$, $0\leq \nu < n$ and $\chi_{n,n}\neq 0$. The computation of the first polynomials of all component sequences is done recursively by the recurrence relations of Proposition \ref{QD-MOPS}.  This computation could also be done  using the general relations of Proposition \ref{prop_nova}, but the ones of Proposition \ref{QD-MOPS} are more efficient.
Studying the properties of the first polynomials, we can obtain some negative answers about the sequences.
In the case some positive affirmations are true for the first polynomials, this give us an indication to try to demonstrate those properties, for all $n$. Furthermore,
in some cases, it is also possible to infer models for the closed formulas to the moments, the structure coefficients or to the recurrence coefficients. Afterwards, we should demonstrate these models by others means.

In the case $\{W_n\}_{n\geq 0}$ is only a MPS, not orthogonal, we start our computations by its structure coefficients $\{\beta_{n}\}_{n \geq 0}$ and $\{\chi_{n,\nu}\} _{0\leq \nu \leq n,\: n \geq 0}$. We will use the recurrence relations of Proposition \ref{prop_nova} for computing recursively the first polynomials of all components sequences. These first polynomials are expanded in the canonical basis.
With respect to the principal components, we can always compute the first moments of their dual sequences and their structure coefficients and wonder if they are orthogonal and symmetric.
Concerning the secondary components, we analyse if they vanish or not, have or not the exact degree, or if they are symmetric or not.
It is worth mentioning here that from the first polynomials, we can also study several other properties.

For each fixed PS $\{W_n\}_{n\geq 0}$, we can test several choices of the parameters $p$, $q$ and $a$ and investigate their effect on the properties of the component sequences. This aspect is particularly interesting when we search for preserving the orthogonality.

Let us summarize step by step the symbolic implementation.

\begin{itemize}
\item  {\bf ORTHOGONAL CASE}: $\{W_n\}_{n\geq 0}$ is a  MOPS

\vspace{0.25cm}

{\bf Step 1 - Starting data }
\vspace{0.1cm}

- the definitions of the recurrence coefficients $\beta_{n}$ and $\gamma_{n+1}$ for all $n\geq 0$ to the MOPS $\{W_n\}_{n\geq 0}$, or their first elements.

- the parameters $p$ and $q$ of $\omega(x)$, the parameter $a$ of $(x-a)$.

- the value of $nmax$.
\vspace{0.1cm}

{\bf Step 2 - For principal components}
\vspace{0.1cm}

- Compute and analyse the closed formulas to the extended coefficients.

\vspace{0.1cm}

{\bf Step 3}

\vspace{0.1cm}

 {\bf - For principal components}

In the case the last two extended coefficients do not vanish, compute the first polynomials,
compute and analyse the first moments of the dual sequences and the first structure coefficients.

\vspace{0.1cm}

{\bf - For secondary components}

Compute and analyse the first polynomials.

\vspace{0.5cm}

\item  {\bf  NON-ORTHOGONAL CASE}: $\{W_n\}_{n\geq 0}$ is a  MPS

\vspace{0.25cm}

{\bf Step 1 - Starting data}

\vspace{0.1cm}

- the definitions of the  structure coefficients $\beta_{n}$ and $\chi_{n,\nu}$ for all $n$, and all $\nu$: $0\leq \nu \leq n$, $n\geq 0$, to the  MPS $\{W_n\}_{n\geq 0}$, or their first elements.

- the parameters $p$, $q$ and $a$.

- $nmax$.

\vspace{0.1cm}

\vspace{0.1cm}

{\bf Step 2}

\vspace{0.1cm}

{\bf  - For principal components}

\vspace{0.1cm}

Compute the first polynomials, compute and analyse the first moments of the dual sequences and the first structure coefficients.

\vspace{0.1cm}

{\bf -  For secondary components}

\vspace{0.1cm}

Compute and analyse the first polynomials.

\vspace{0.25cm}

\item
For each sequence $\{W_n\}_{n\geq 0}$, test several  values of parameters $p$, $q$ and $a$.
\end{itemize}

\section{Study cases}\label{results}

In this section, we present new results obtained with the symbolic implementation for several cases of orthogonal sequences.

\subsection{A symmetric case}
Let $\left\{W_n\right\}_{n\geq 0}\:$ be a symmetric, semi-classical of class 1 MOPS \cite{AlP:96,C:78,PZ:13}  with coefficients
\begin{eqnarray}
&& \beta_n=0\ ,\ n\geq 0\ , \notag\\
&&\gamma_{2n+1}=\displaystyle\frac{(\beta+n+1)(\alpha+\beta+n+1)}{(\alpha+\beta+2n+1)(\alpha+\beta+2n+2)} \ ,\ n\geq 0\ ,\notag\\
&& \gamma_{2n+2}=\displaystyle\frac{(n+1)(\alpha+n+1)}{(\alpha+\beta+2n+2)(\alpha+\beta+2n+3)}\ ,\ n\geq 0\ , \notag
\end{eqnarray}
and regularity conditions $\alpha\neq-(n+1)$, $\beta\neq-(n+1)$, $\alpha+\beta\neq-(n+1)$, $n\geq 0$. From symmetrie, we know that $(w_0)_{2n+1}=0$, $n\geq 0$.
The extended coefficients of $\{P_n\}_{n\geq 0}$ and $\{R_n\}_{n\geq 0}$ are
\begin{eqnarray}
\beta_0^P & = & a p+q+\frac{\beta +1}{\alpha +\beta +2} \notag \\
\beta_{n+1}^P & = & a^2+a p-a (a+p)+q+\frac{(n+1) (n+\alpha +1)}{(2 n+\alpha +\beta +2)
(2 n+\alpha +\beta +3)}\notag\\
&& +\frac{(n+\beta +2)
(n+\alpha +\beta +2)}{(2 n+\alpha +\beta +3) (2
n+\alpha+\beta+4)}\notag\\
\gamma_{n+1}^P & = & \frac{(n+1) (n+1)
(n+\alpha +1) (n+\alpha +1)}{(2 n+\alpha +\beta +1) (2 n+\alpha
+\beta +2)^2 (2 n+\alpha +\beta +3)} \notag \\
\varrho_{n+1}^P & = & p \quad , \quad \rho_{n+1}^P=p\frac{(n+1) (n+\alpha +1)}{(2 n+\alpha +\beta +2) (2 n+\alpha +\beta +3)} \notag
\end{eqnarray}
\begin{eqnarray}
\beta_0^R & = & \frac{-p^2 (\alpha +\beta +3)+q (\alpha +\beta +3)+\beta +2}{\alpha +\beta +3} \notag \\
\beta_{n+1}^R & = & \frac{n^2 (4 q+2)+2 n (2 q+1) (\alpha +\beta +4)+
\alpha  \beta +4 \alpha +\beta^2+5 \beta +8}
{(2n+\alpha +\beta +3) (2 n+\alpha +\beta +5)} \notag \\
&& +q\frac{\left(\alpha^2+2 \alpha  (\beta +4)+\beta^2+8 \beta
+15\right)}{(2
n+\alpha +\beta +3) (2 n+\alpha +\beta +5)}\notag\\
\gamma_{n+1}^R & = & \frac{(n+\beta +2) (n+\alpha +\beta +2) (n+1) (n+\alpha +1)}
{(2 n+\alpha +\beta +2) (2 n+\alpha +\beta +3)^2 (2 n+\alpha
+\beta +4)} \notag\\
\varrho_{n+1}^R & = & p \quad,\quad \rho_{n+1}^R=p\frac{(n+\beta +2) (n+\alpha +\beta +2)}{(2 n+\alpha +\beta +3) (2 n+\alpha +\beta +4)} \notag
\end{eqnarray}


Choosing the parameter  $\:p=0\:$, we obtain orthogonality of the principal component sequences with the following recurrence coefficients obtained from the preceding formulas after some simplifications
\begin{eqnarray}
\beta_{0}^{P}& = & 1+q-\displaystyle\frac{\alpha+1}{\alpha+\beta+2} \notag \\
\beta_{n+1}^{P} & = & 1+q+\displaystyle\frac{(n+1)\left(n+\alpha+1\right)}{\alpha+\beta+2n+2}
-\displaystyle\frac{(n+2)
\left(n+2+\alpha\right)}{\alpha+\beta+2n+4}  \notag \\
\gamma_{n+1}^{P} & = & \displaystyle\frac{(n+1)(\alpha+n+1)(\beta+n+1)(\alpha+\beta+n+1)}
{(\alpha+\beta+2n+1)(\alpha+\beta+2n+2)^2(\alpha+\beta+2n+3)}\notag \\
\beta_{0}^{R} & = &1+q-\displaystyle\frac{\alpha+1}{\alpha+\beta+3}  \notag \\
\beta_{n+1}^{R} & = &1+q+\displaystyle\frac{(n+1)\left(\alpha+n+1\right)}{\alpha+\beta+2n+3}
-\displaystyle\frac{(n+2)\left(\alpha+n+2\right)}{\alpha+\beta+2n+5}  \notag \\
\gamma_{n+1}^{R} & = & \displaystyle\frac{(n+1)(\alpha+n+1)(\beta+n+2)(\alpha+\beta+n+2)}
{(\alpha+\beta+2n+2)(\alpha+\beta+2n+3)^2(\alpha+\beta+2n+4)}  \notag
\end{eqnarray}
We remark that the recurrence coefficients of $\{R_n\}_{n\geq 0}$ can be obtained from the ones of $\{P_n\}_{n\geq 0}$ replacing $\beta$ by $\beta+1$.

For $n=0,1,\ldots,nmax$, we obtain
$a_n(x)=0$, $b_n(x)=a\, R_n(x)$,
which leads us to infer that the sequence
$\:\left\{b_n\right\}_{n\geq 0}\:$ is orthogonal and the sequence $\:\left\{a_n\right\}_{n\geq 0}\:$ is void.
It is possible to prove this for all $n\geq 0$, using the relations (\ref{sist6.1}) and (\ref{sist6.2}) of Proposition \ref{QD-MOPS-anbn}.
In fact, in this case, the GQD is
\begin{equation} \label{ex1W2nW2n1}
W_{2n}(x)=P_n(x^2-q)\quad,\quad W_{2n+1}(x)=xR_n(x^2-q)\ .
\end{equation}
From this and Proposition 4.1 in \cite[p.1326]{A:2010} it is
possible to deduce that
\begin{equation} \label{ex1w02nu02n}
(u_0)_{n}= \sum_{k=0}^{n} \binom{n}{k} \, q^{n-k}(w_0)_{2k}\ .
\end{equation}


With respect to the canonical forms $u_{0}$ and $v_{0}$ of the principal components, their first moments for $nmax=3$ are
\medskip
\begin{eqnarray}
(u_0)_{n=0,1,2,3} & = &  \left\{ 1\,,\, q+\displaystyle\frac{\beta +1}{\alpha +\beta +2}\,,\, \right.  \notag\\
&& \left. \ \ \ q^2+ \frac{2q(\beta +1)}{\alpha +\beta +2}+
\frac{(\beta +1) (\beta +2)}{(\alpha +\beta +2) (\alpha +\beta +3)} \,,\,\right.\notag\\
&& \ \ \ \left. q^3 + \displaystyle\frac{3q^2 (\beta +1)}{\alpha +\beta +2 } + \frac{3q (\beta +1)(\beta +2)}{(\alpha +\beta +2 )(\alpha +\beta +3)}+\right. \notag\\
&& \ \ \ \ \left. \frac{(\beta +1)(\beta +2)(\beta +3)}{(\alpha +\beta +2 )(\alpha +\beta +3)(\alpha +\beta +3)} \right\} \ ,\notag\\
(v_0)_{n=0,1,2,3}& = & \left\{1\,,\, q+\displaystyle\frac{\beta +2}{\alpha +\beta +3}  \,,\, \right.\notag\\
&& q^2+ \frac{2q(\beta +2)}{\alpha +\beta +3}+
\frac{(\beta +2) (\beta +3)}{(\alpha +\beta +3) (\alpha +\beta +4)} \,,\,   \notag\\
&& \ \ \ \left. q^3 + \displaystyle\frac{3q^2 (\beta +2)}{\alpha +\beta +3 } + \frac{3q (\beta +2)(\beta +3)}{(\alpha +\beta +3 )(\alpha +\beta +4)}+ \right. \notag\\
&& \ \ \ \left.  \frac{(\beta +2)(\beta +3)(\beta +4)}{(\alpha +\beta +3 )(\alpha +\beta +4)(\alpha +\beta +5)} \right\}  \ .\notag
\end{eqnarray}

Analysing the moments for higher values of $nmax$, allow us to infer the next models for the closed formulas valid for all $n\geq 0$,
supposing  that $ \prod_{i=1}^{0}\left( .\right) = 1$.
\begin{eqnarray}
\left(u_{0}\right)_{n} & \leftarrow & \sum_{k=0}^{n} \binom{n}{k} \, q^{n-k} \frac{\displaystyle \prod_{i=1}^{k}\left( \beta +i\right)}{\displaystyle \prod_{i=1}^{k}\left( \alpha + \beta +1+ i \right)}\ ,\label{u0nSC}
\end{eqnarray}
\begin{eqnarray}
\left(v_{0}\right)_{n} & \leftarrow & \sum_{k=0}^{n} \binom{n}{k} \, q^{n-k} \frac{\displaystyle \prod_{i=1}^{k}\left( \beta +1+i\right)}{\displaystyle \prod_{i=1}^{k}\left( \alpha + \beta +2+ i \right)}\ .\notag
\end{eqnarray}
From (\ref{ex1w02nu02n}) and (\ref{u0nSC}), we can also identify the
following model
\begin{eqnarray}
(w_0)_{2k} & \leftarrow & \frac{\displaystyle \prod_{i=1}^{k}\left( \beta +i\right)}{\displaystyle \prod_{i=1}^{k}\left( \alpha + \beta +1+ i \right)}\ .\notag
\end{eqnarray}
Afterwards, we should demonstrate these formulas by other means.

\subsection{A non-symmetric case}

Let $\left\{W_n\right\}_{n\geq 0}\:$ be MOPS with constant recurrence coefficients, also considered in \cite{T:2012} in the context of cubic decomposition
$$\beta_n=\beta \quad; \quad  \gamma_{n+1}=\gamma \neq 0\ ,\quad  n\geq 0\ .$$
It is worth noting that this MOPS can be obtained by applying a shifting transformation to the monic Chebyshev sequence of second kind $\{U_{n}\}_{ n \geq 0}$ whose recurrence coefficients are $\:\beta^{U}_{n}=0\:$ and $\:\gamma^{U}_{n+1}=\frac{1}{4} \,,\: n\geq 0\:$.
From identities (\ref{shifted-rec-coeff}), we obtain
$$W_{n}\left(x\right) = A^{-n} U_{n}\left( Ax+B  \right)\quad , \quad Ax+B = \pm\frac{1}{2\sqrt{\gamma}}\left( x-\beta \right).$$
Since the Chebyshev sequence of second kind belongs to the Jacobi class of classical polynomials, and the shifting transformation preserves the class \cite{P:94}, we know that $\left\{W_n\right\}_{n\geq 0}\:$ remains in the Jacobi class.
The extended coefficients of $\{P_n\}_{n\geq 0}$ and $\{R_n\}_{n\geq 0}$ are

\vspace{0.35cm}

\begin{tabular}{ll}
$\beta_0^P=q+\gamma-\beta^2+a (p+2 \beta)$   & $,\quad \beta_0^R=q+2\gamma-3\beta\left(p+\beta\right)-p^2$\ , \notag\\
$\beta_{n+1}^P=q+2 \gamma+\beta (p+\beta)$  & $,\quad \beta_{n+1}^R=q+2\gamma+\beta(p+\beta)$\ , \notag \\
$\gamma_{n+1}^P=\gamma^2$  & $,\quad \gamma_{n+1}^R=\gamma^2$\ , \notag \\
$\varrho_{n+1}^P=p+2\beta$   & $,\quad \varrho_{n+1}^R=p+2\beta$\ , \notag \\
$\rho_{n+1}^P=\gamma\left(p+2\beta\right)$   & $,\quad \rho_{n+1}^R=\gamma(p+2\beta)$\ . \notag
\end{tabular}


\vspace{0.35cm}

Taking the parameter  $\:p=-2\beta\:$, the principal components are orthogonal with the following recurrence coefficients

\vspace{0.35cm}
\begin{tabular}{ll}
$\beta_{n}^R=q+2\gamma-\beta^2$ & $,\ \gamma_{n+1}^R=\gamma^2\ ,\quad n\geq 0\ .$\\
$\beta_0^P=\beta_0^R-\gamma\ ,\ \beta_{n+1}^P=\beta_{n+1}^R\ $ & $,\ \gamma_{n+1}^P=\gamma_{n+1}^R\ , \ n\geq 0\ .$\\
\end{tabular}
\vspace{0.35cm}

These identities show that if all parameters are real, then the principal components are both positive definite, because $\beta^R_{n},\ \beta^P_{n}\in {\bkR}$ and
$\gamma^R_{n+1}>0$, $n\geq 0$.
Furthermore, we see that
$\{P_n\}_{n\geq 0}$ is a co-recursive of $\{R_n\}_{n\geq 0}$ with perturbation $-\gamma$, i.e.,
\begin{equation} \label{xxx}
P_n(x)=R_n(-\gamma;x)\ ,\quad n\geq 0.
\end{equation}
Also, we can state that $\{R_n\}_{n\geq 0}$ is a shifting of $\{W_n\}_{n\geq 0}$ as follows
$$
R_{n}\left(x\right) = A^{-n} W_{n}\left( Ax+B  \right)\quad, \quad A=\pm\sqrt{\frac{1}{\gamma}}\quad, \quad B=-A(q+2\gamma-\beta^2)+\beta \ .$$
For $n=0,1,\ldots,nmax$, we obtain that
$a_n(x)=0$, $b_n(x)=(a-\beta)R_n(x)$,
which allows us to infer that the sequence $\{b_{n}(x)\}_{n\geq 0}$ is orthogonal and the sequence $\{a_{n}(x)\}_{n\geq 0}$ is void, which
can be proven by Proposition \ref{QD-MOPS-anbn}. In fact, in this case, we have (\ref{xxx}) and thus the GQD of $\left\{W_n\right\}_{n\geq 0}\:$ is given by
\cite{ATese:2004,A:2010}
$$W_{2n}(x)=R_n\big(-\gamma;\omega (x)\big)\quad ,\quad
W_{2n+1}(x)=\big(x-\beta\big)R_n\left(\omega (x)\right) \;,\:n\geq 0\,.\:$$

Regarding the canonical forms $u_{0}$ and $v_{0}$ of the principal components, we obtained for $nmax=3$ the next first moments
\begin{eqnarray}
(u_0)_{n=0,1,2,3} & = & \left\{\ 1\ ,\ \gamma -\beta^2+q\ ,\ \ 2 \gamma^2+\beta^4+q^2-2 \beta^2 (\gamma +q)+2\gamma  q\ ,\right.  \notag\\
&& \ \ \left. 5 \gamma^3-\beta^6+q^3-3 \beta^2 \left(2 \gamma^2+q^2+2 \gamma  q\right)+ \right.\notag \\
&& \ \ \left.  3 \gamma  q^2+6 \gamma^2 q+3 \beta^4 (\gamma+q) \ \right\}.\notag \\
(v_0)_{n=0,1,2,3} & = & \left\{\ 1\ ,\  2 \gamma -\beta^2+q\ ,\ \ 5 \gamma^2+\beta^4+q^2-2 \beta^2 (2 \gamma +q)+4 \gamma  q\ ,\right.\notag\\
&& \ \ \left. \left(2 \gamma -\beta^2+q\right) \left(7 \gamma^2+\beta ^4+q^2-2 \beta^2 (2 \gamma +q)+4 \gamma  q \right)\  \right\}.\notag
\end{eqnarray}

\subsection{A collection of orthogonal cases}

In this section we present explicit results obtained by applying the symbolic implementation described in Section \ref{sectionSIGQD} to a collection of well known
orthogonal sequences. We consider all classical families of Hermite, Laguerre, Bessel and Jacobi \cite{P:94}, the generalized Hermite sequence \cite{C:78}, which is
an example of symmetric semi-classical sequence of class one, and the classical discrete Charlier sequence \cite{C:78}. We give  results for particular cases of Gegenbauer, namely for Legendre and all Chebyshev sequences, in order to be able to identify some
existing relationships between the original sequences and their components.

From the recurrence coefficients of each MOPS $\{W_n\}_{n\geq0}$ considered, we have computed the extended coefficients of the principal components $\{P_n\}_{n\geq0}$ and $\{R_n\}_{n\geq0}$ with parameters $p$, $q$ and $a$ arbitrary.
After that, we have fixed the values of the parameters considering the cases for which they are all equal to zero ($a=p=q=0$), two of them are zero ($p=q=0$, $a=p=0$ or $a=q=0$) and only one of them is zero ($a=0$, $p=0$ or $q=0$), resulting in a total of seven cases to study for each given sequence.
In each case we signalize the orthogonality or non-orthogonality of the components and
we point out vanishing secondary sequences. For Chebyshev polynomials, we have taken others particular values of parameters that lead to simplifications and we mention some relationships between the sequences with respect to shifting, association and  perturbation transformations. All parameters inherent to each sequence were left free.
We remark that in the cases of Hermite, generalized Hermite and Gegenbauer (including Legendre and all Chebyshev) sequences the parameter $p=0$ leads to
orthogonality of the principal components.

\section{Conclusions}

The symbolic implementation described along this work allows an automatic enquiry regarding the properties of orthogonality and symmetry of the component sequences
emerging from a GQD of a given MPS.
We can also test several choices of the parameters and conclude about their effect on the properties of the components.
In the case, the original sequence is orthogonal, we explicit the closed formulas of the extended coefficients of the principal component sequences.
We recall that the extended  coefficients are reduced to the recurrence coefficients when the principal components are orthogonal,
furthermore they allow an efficient way of computation of the first polynomials of all component sequences.

From the collection of results obtained, we point out the interest of the general decomposition, even in some symmetric examples.
In fact, for these examples, when $p=0$ the orthogonality of the original sequence is preserved in such a way the principal components are also orthogonal.

In conclusion, we consider that this implementation is an efficient method for computing the GQD, testing some properties and direct the theoretical study about this kind of decomposition.

\section*{Acknowledgements}

T. A. Mesquita and Z. da Rocha were partially supported by
\newline CMUP (UID/MAT/00144/2013), which is funded by FCT (Portugal) with national (MEC) and European
structural funds (FEDER), under the partnership agreement PT2020.

\noindent
\^A. Macedo was financed by Portuguese Funds through FCT - Funda\c{c}{\~a}o para a Ci{\^e}ncia e a Tecnologia, within the Project UID/MAT/00013/2013.

\newpage
%
%
%

\begin{center}
\begin{tabular}[c]{p{3.0cm}||p{7.5cm}}
 &  {\bf  CHARLIER } \\[.5em]
\hline \hline
$\beta_{n}\ \ | \ \ \gamma_{n+1}$ &  $n+\alpha\ \ |\ \  (n+1)\alpha$\\[.5em]
\hline \hline
$\beta_{n}^P$ & $\beta_0^P=a(p+2\alpha+1)-\alpha^2$     \\[.5em]
&
$\beta_{n+1}^P=4 n^2+2 n (p+4 \alpha +4)+p (\alpha +2)+\alpha ^2+9 \alpha +4$    \\[.8em]
\hline
$\gamma_{n+1}^P\ \ |\ \ \varrho_{n+1}^P$ & $2(n+1)(2 n+1)\alpha ^2\ \ |\ \ 4n+p+2 \alpha +5$
 \\[.5em]
\hline
 $\rho_{n+1}^P$ & $2(n+1)\alpha(4 n+p+2\alpha+3)$  \\[.5em]
\hline \hline $\beta_{n}^R$ &
$\beta_0^R=-p^2-3p(\alpha +1)-3\alpha^2-3\alpha-2$ \\[.5em]
& $\beta_{n+1}^R=4 n^2+2 n (p+4 \alpha +6)+p (\alpha +3)+\alpha^2+13\alpha +9$   \\[.8em]
\hline
$\gamma_{n+1}^R\ \ |\ \ \varrho_{n+1}^R$ & $2 (n+1) (2 n+3) \alpha ^2\ \ |\ \ 4 n+p+2 \alpha +7$  \\[.5em]
\hline $\rho_{n+1}^R$ & $(2 n+3) \alpha  (4 n+p+2 \alpha +5)$   \\[.5em]
\hline \hline
all cases &$P_n$, $R_n$, $a_n$, $b_n$ nonorthogonal \\
\hline
\end{tabular}
\end{center}

\begin{center}
\begin{tabular}[c]{p{3.0cm}||p{7.5cm}}
 & {\bf GENERALIZED HERMITE}\\[.5em]
\hline \hline
$\beta_{n}\ \ | \ \ \gamma_{n+1}$
& $\beta_n=0\ \ | \ \ \gamma_{n+1}=\frac{1}{2}\left(n + 1 + \mu\left(1 +(-1)^n\right)\right)$\\[.7em]
\hline \hline
$\beta_{n}^P$
&$\beta_0^P=\frac{1}{2}+a p+q+\mu\ ,\ \beta_{n+1}^P=\frac{5}{2}+2 n+q+\mu$
\\[.5em]
\hline $\gamma_{n+1}^P\ \ | \ \ \varrho_{n+1}^P  \ \ | \ \  \rho_{n+1}^P$ &$\frac{1}{2}(n+1)(2 n+2\mu +1)\ \ | \ \ p\ \ |\ \ p \left(n+1\right)$
\\[.5em]
\hline \hline $\beta_{n}^R$ &
$\beta_0^R=\frac{3}{2}-p^2+q+\mu\ ,\ \beta_{n+1}^R=\frac{7}{2}+2 n+q+\mu$
\\[.5em]
\hline $\gamma_{n+1}^R\ \ | \ \ \varrho_{n+1}^R\ \ |\ \ \rho_{n+1}^R$ &$\frac{1}{2} \left(n+1\right) \left(2n+2\mu
+3\right)\ \ | \ \ p\ \ | \ \ \frac{1}{2} p \left(2 n+2\mu +3\right)$
\\[.5em]
\hline\hline
 $p=0$ &$P_n$, $R_n$, $b_n=aR_n$ orthogonal;  $a_n=0$ \\
\hline\hline
$a=p=q=0$
&$P_n$, $R_n$ orthogonal; $a_n=0$, $b_n=0$ \\
\hline $p=q=0$
&$P_n$, $R_n$, $b_n=aR_n$ orthogonal;  $a_n=0$  \\
\hline $a=p=0$ &$P_n$, $R_n$ orthogonal;  $a_n=0$ and $b_n=0$
\\
\hline $a=q=0$ &$P_n$, $R_n$, $a_n$, $b_n$ nonorthogonal
\\
\hline $a=0$ &$P_n$, $R_n$, $a_n$, $b_n$ nonorthogonal
\\
\hline $q=0$ &$P_n$, $R_n$, $a_n$, $b_n$ nonorthogonal
\\
\hline
\end{tabular}
\end{center}

\begin{center}
\begin{tabular}[c]{p{3.0cm}||p{7.5cm}}
 &  {\bf HERMITE } \\[.5em]
\hline \hline
$\beta_{n}\ \ | \ \ \gamma_{n+1}$ & $0 \ \ | \ \ \frac 12(n+1)$ \\[.5em]
\hline \hline $\beta_{n}^P$ & $\beta_0^P=\frac 12+ap+q\ , \ \beta_{n+1}^P=\frac 12+2(n+1)+q$  \\[.5em]
\hline $\gamma_{n+1}^P\ \ | \ \ \varrho_{n+1}^P\ \ | \ \ \rho_{n+1}^P$ &
$\frac 12(n+1)(2n+1)\ \ | \ \ p\ \ | \ \ (n+1)p$   \\[.5em]
\hline \hline $\beta_{n}^R$ & $\beta_0^R=\frac 32-p^2+q\ ,\  \beta_{n+1}^R=\frac 32+2(n+1)+q$  \\[.5em]
\hline $\gamma_{n+1}^R\ \ | \ \ \varrho_{n+1}^R  \ \ | \ \ \rho_{n+1}^R$ &
$\frac 12(n+1)(2n+3)  \ \ | \ \  p \ \ | \ \  \left(n+\frac 32\right)p $ \\[.5em]
\hline \hline
$p=0$ & $P_n$, $R_n$, $b_n$,  $a_n$ orthogonal \\
\hline
$a=p=q=0$ & $P_n$, $R_n$ orthogonal ; $a_n=0$, $b_n=0$ \\
\hline
$p=q=0$ & $P_n$, $R_n$, $b_n=aR_n$ orthogonal ; $a_n=0$  \\
\hline
$a=p=0$ & $P_n$, $R_n$ orthogonal ; $a_n=0$, $b_n=0$ \\
\hline
$a=q=0$ &$P_n$, $R_n$, $a_n$, $b_n$ nonorthogonal  \\
\hline
$a=0$ & $P_n$, $R_n$, $a_n$, $b_n$ nonorthogonal \\
\hline
$q=0$ & $P_n$, $R_n$, $a_n$, $b_n$ nonorthogonal \\
\hline
\end{tabular}
\end{center}
%

\vfill
\begin{center}
\begin{tabular}[c]{p{3.0cm}||p{7.5cm}}
 & {\bf  LAGUERRE }\\[.5em]
\hline \hline
$\beta_{n} \ \ | \ \ \gamma_{n+1}$ & $2n+\alpha+1\ \ | \ \ (n+1)(n+\alpha +1)$ \\[.5em]
\hline \hline $\beta_{n}^P$ & $\beta_0^P=-2+q-3\alpha-\alpha^2+a(4+p+2\alpha)$ \\[.5em]
&$\beta_{n+1}^P=2+p+q+\alpha(3+p+\alpha)+$\\
 & $\ \ \ \ \ \ \ \ \ \  4(n+1)(9+6n+p+3\alpha)$    \\[.8em]
\hline $\gamma_{n+1}^P$ & $ 2(n+1) (2n+1)
(2 n+\alpha+1)(2 n+\alpha+2)$  \\[.5em]
\hline
$\varrho_{n+1}^P$ & $8 n+p+2(\alpha+ 6)$    \\[.5em]
\hline $\rho_{n+1}^P$ & $2(n+1) (2 n+\alpha+2) (8+8
n+p+2\alpha)$ \\[.5em]
\hline \hline $\beta_{n}^R$ &
$\beta_0^R=-p^2+q-3 p (3+\alpha)-3 (6+5 \alpha+\alpha^2)$  \\[.5em]
& $\beta_{n+1}^R=14+q +3p+\alpha(\alpha+9+p)+$\\[.8em]
& $\ \ \ \ \ \ \ \ \ \ \ 4(n+1) (6n +15 + p + 3 \alpha)$\\[.8em]
\hline $\gamma_{n+1}^R$ & $2(1+n)(3+2n)(2+2 n+\alpha)(3+2 n+\alpha)$\\[.5em]
\hline $\varrho_{n+1}^R$ & $8n+p+2(8+\alpha)$
 \\[.5em]
\hline $\rho_{n+1}^R$ & $(3+2 n) (3+2n+\alpha) (8 n+p+2 (6+\alpha))$    \\[.5em]
\hline \hline
all cases &$P_n$, $R_n$, $a_n$, $b_n$ nonorthogonal \\
\hline
\end{tabular}
\end{center}
\vfill
\begin{center}
\begin{tabular}[c]{p{1.5cm}||p{9.0cm}}
 &  {\bf  BESSEL } \\[.5em]
\hline \hline
$\beta_{n}$ &  $\beta_0=-\frac 1\alpha$    \\[.5em]
&   $\beta_{n+1}=\frac{1-\alpha}{(n+\alpha)(n+\alpha+1)}$    \\[.8em]
$\gamma_{n+1}$ &    $-\frac{(n+1)(n+2\alpha-1)}{(2n+2\alpha-1)(n+\alpha)^2(2n+2\alpha+1)}$  \\[.5em]
\hline \hline & \\
$\beta_{n}^P$ & $\beta_0^P=\frac{-2+q+3 q \alpha+2 q \alpha^2+a (1+2 \alpha) (-2+p+p \alpha)}{1+3 \alpha+2 \alpha^2}$     \\[.5em]
&
$\beta_{n+1}^P=\Big(64 n^4 q+64 n^3 q (2 \alpha +3)+$\\
& $\ \ \ \ \ \  \ \ \  \ \ \ 4 n^2 \left(-4p (\alpha -1)+q \left(24 \alpha ^2+72 \alpha
   +49\right)-2\right)+$\\
   & $\ \ \ \ \ \  \ \ \  \ \ \ 2 n (2 \alpha +3) \left(-4 p (\alpha -1)+q \left(8 \alpha ^2+24 \alpha
   +13\right)-2\right)+$\\
   & $\ \ \ \ \ \  \ \ \  \ \ \ p \left(-4 \alpha ^3-8 \alpha ^2+7 \alpha +5\right)+$\\
   & $\ \ \ \ \ \  \ \ \  \ \ \ 4 q \alpha ^4+24 q \alpha ^3+49
   q \alpha ^2+39 q \alpha +10 q+4 \alpha ^2-18 \alpha +2\Big)/$\\
   & $\ \ \ \ \ \  \ \ \  \ \ \ \Big((2 n+\alpha +1) (2 n+\alpha +2) (4 n+2 \alpha
   +1) (4 n+2 \alpha +5)\Big)$    \\[.8em]
\hline
& \\
$\gamma_{n+1}^P$ &$\frac{4 (n+1) (2 n+1) (n+\alpha )(2 n+2 \alpha -1)}{(4 n+2 \alpha -1) (4 n+2 \alpha +1)^2 (4n+2\alpha
   +3) \left(4 n^2+4 n \alpha +2 n+\alpha ^2+\alpha \right)^2}$
 \\[.5em]
\hline
& \\
$\varrho_{n+1}^P$ &$\frac{2-2 \alpha+p (3+4 n^2+4 \alpha+\alpha^2+4 n (2+\alpha))}{(1+2 n+\alpha) (3+2 n+\alpha)}$   \\[.5em]
\hline
& \\
$\rho_{n+1}^P$ &
$-\frac{4 (1+n) (n+\alpha) (2+4 n^2 p+2 (-1+p) \alpha+p \alpha^2+4 n p (1+\alpha))}
{(2 n+\alpha) (1+2 n+\alpha)^2 (2+2 n+\alpha) (1+4 n+2 \alpha) (3+4 n+2 \alpha)}$  \\[.5em]
\hline \hline
& \\
$\beta_{n}^R$ &
$\beta_0^R=\frac{-6+p (9+6 \alpha)-p^2 (6+7 \alpha+2 \alpha^2)+q (6+7 \alpha+2 \alpha^2)}{6+7 \alpha+2 \alpha^2}$ \\[.5em]
&
$\beta_{n+1}^R=\Big(-6+126 q+64 n^4 q-22 \alpha+225 q \alpha+4
\alpha^2+145 q \alpha^2+40 q \alpha^3+4 q \alpha^4+64 n^3 q (5+2
\alpha)-p (-21+\alpha+16 \alpha^2+4 \alpha^3)+2 n (5+2 \alpha)
(-2-4 p (-1+\alpha)+q (45+40 \alpha+8 \alpha^2))+$\\
&$4 n^2
(-2-4 p (-1+\alpha)+q (145+120 \alpha+24 \alpha^2))\Big)/$\\
& $\Big((2+2 n+\alpha) (3+2 n+\alpha) (3+4 n+2 \alpha) (7+4 n+2 \alpha)\Big)$   \\[.8em]
\hline
& \\
$\gamma_{n+1}^R$ &$\frac{4 (1+n) (3+2 n) (n+\alpha)
(1+2 n+2 \alpha)}
{(1+4 n+2 \alpha) (3+4 n+2 \alpha)^2 (5+4 n+2 \alpha) (2+4 n^2+3 \alpha+\alpha^2+n (6+4 \alpha))^2}$  \\[.5em]
\hline $\varrho_{n+1}^R$ & $\frac{p \left(4 n^2+4 n (\alpha +3)+\alpha ^2+6 \alpha +8\right)-2 \alpha +2}{(2 n+\alpha +2) (2 n+\alpha
   +4)}$   \\[.5em]
\hline $\rho_{n+1}^R$ &$-\frac{(3+2 n) (1+2 n+2 \alpha)
(2-2 \alpha+p (3+4 n^2+4 \alpha+\alpha^2+4 n (2+\alpha)))}
{(1+2 n+\alpha) (2+2 n+\alpha)^2 (3+2 n+\alpha) (3+4 n+2 \alpha) (5+4 n+2 \alpha)}$   \\[.5em]
\hline \hline
all cases &$P_n$, $R_n$, $a_n$, $b_n$ nonorthogonal \\
\hline
\end{tabular}
\end{center}

\begin{center}
\begin{tabular}[c]{p{1.25cm}||p{9.75cm}}
 & {\bf JACOBI}\\
\hline \hline
$\beta_{n}$ & $\beta_0=\frac{\alpha - \beta}{\alpha + \beta + 2}\ \ |\ \
\beta_{n+1}=\frac{\alpha^2 - \beta^2}{(2n + \alpha + \beta + 1)(2n + \alpha + \beta + 2)^2(2n + \alpha + \beta + 3)}$ \\[.7em]
$\gamma_{n+1}$ & $\frac{4(n + 1)(n +\alpha +\beta+1)(n+\alpha+1)(n + \beta + 1)}{(2n +\alpha+\beta+ 1)(2n + \alpha +\beta + 2)^2(2n +\alpha+\beta+3)}$
\\[.5em]
\hline \hline $\beta_{n}^P$ & $\beta_0^P=\omega(a)+\frac{4(1+\alpha) (1+\beta)}{(2+\alpha+\beta)^2 (3+\alpha+\beta)}
-\left(a+\frac{-\alpha+\beta}{2+\alpha+\beta}\right)$\\
& $\left(a+\frac{-\alpha^2+\beta^2}{(3+\alpha+\beta)
(4+\alpha+\beta)^2 (5+\alpha+\beta)}\right)$
\\[.5em]
&$\beta_{n+1}^P=\omega(a)+\frac{8 (1+n) (2+2 n+\alpha) (2+2n+\beta) (2+2 n+\alpha+\beta)} {(3+4 n+\alpha+\beta)
(4+4n+\alpha+\beta)^2 (5+4 n+\alpha+\beta)}+$\\
& $\frac{4 (3+2 n)(3+2n+\alpha) (3+2 n+\beta) (3+2 n+\alpha+\beta)}
{(5+4n+\alpha+\beta)(6+4 n+\alpha+\beta)^2(7+4 n+\alpha+\beta)}
-$\\
& $\left(a+p+\frac{(\alpha-\beta) (\alpha+\beta)}{(5+4n+\alpha+\beta)(6+4 n+\alpha+\beta)^2(7+4 n+\alpha+\beta)}\right)\times$\\
&
$\left(a+\frac{-\alpha^2+\beta^2}{(5+4 n+\alpha+\beta)(6+4n+\alpha+\beta)^2 (7+4 n+\alpha+\beta)}\right)$
\\[.8em]
\hline $\gamma_{n+1}^P$ & $\frac{32 (1+n) (1+2 n) (1+2 n+\alpha)
(2+2 n+\alpha) (1+2 n+\beta) (2+2 n+\beta) (1+2 n+\alpha+\beta)
(2+2 n+\alpha+\beta)} {(1+4 n+\alpha+\beta) (2+4 n+\alpha+\beta)^2
(3+4 n+\alpha+\beta)^2 (4+4 n+\alpha+\beta)^2 (5+4n+\alpha+\beta)}$
\\[.5em]
\hline $\varrho_{n+1}^P$ & $p+\frac{(\alpha-\beta)
(\alpha+\beta)}{(5+4 n+\alpha+\beta) (6+4 n+\alpha+\beta)^2 (7+4n+\alpha+\beta)} +$\\
& $\frac{(\alpha-\beta) (\alpha+\beta)}
{(7+4n+\alpha+\beta) (8+4 n+\alpha+\beta)^2 (9+4n+\alpha+\beta)}$
\\[.5em]
\hline $\rho_{n+1}^P$ &$8 (1+n) (2+2n+\alpha) (2+2n+\beta)(2+2n+\alpha+\beta)\times$\\
& $\left(\frac{p+\frac{(\alpha-\beta)
(\alpha+\beta)} {(3+4 n+\alpha+\beta)(4+4n+\alpha+\beta)^2(5+4n+\alpha+\beta)} +\frac{(\alpha-\beta) (\alpha+\beta)}
{(5+4n+\alpha+\beta) (6+4 n+\alpha+\beta)^2 (7+4 n+\alpha+\beta)}}
{(3+4 n+\alpha+\beta) (4+4n+\alpha+\beta)^2 (5+4n+\alpha+\beta)}\right)$
\\[.5em]
\hline \hline $\beta_{n}^R$ & $\beta_0^R=\omega(a)+\frac{4(1+\alpha) (1+\beta)}{(2+\alpha+\beta)^2 (3+\alpha+\beta)}
+\frac{8 (2+\alpha) (2+\beta) (2+\alpha+\beta)}{(3+\alpha+\beta)(4+\alpha+\beta)^2 (5+\alpha+\beta)}-$\\
& $\left(p+\frac{\alpha-\beta}{2+\alpha+\beta}+\frac{(\alpha-\beta)
(\alpha+\beta)}{(3+\alpha+\beta) (4+\alpha+\beta)^2
(5+\alpha+\beta)} \right)$\\
& $\left(a+p+\frac{(\alpha-\beta)(\alpha+\beta)}{(5+\alpha+\beta) (6+\alpha+\beta)^2
(7+\alpha+\beta)}\right)
-$\\
& $\left(a+\frac{-\alpha+\beta}{2+\alpha+\beta}\right)\left(a+\frac{-\alpha^2+\beta^2}{(3+\alpha+\beta)
(4+\alpha+\beta)^2 (5+\alpha+\beta)} \right)$
\\[.5em]
& $\beta_{n+1}^R=\omega(a)+\frac{4 (3+2 n) (3+2 n+\alpha) (3+2
n+\beta) (3+2 n+\alpha+\beta)}{(5+4 n+\alpha+\beta) (6+4
n+\alpha+\beta)^2 (7+4 n+\alpha+\beta)} +$\\
& $\frac{8 (2+n) (4+2
n+\alpha) (4+2 n+\beta) (4+2 n+\alpha+\beta)}{(7+4 n+\alpha+\beta)
(8+4 n+\alpha+\beta)^2 (9+4 n+\alpha+\beta)}
-$\\
& $\left(a+p+\frac{(\alpha-\beta) (\alpha+\beta)}{(7+4
n+\alpha+\beta) (8+4 n+\alpha+\beta)^2 (9+4
n+\alpha+\beta)}\right)\times$\\
& $\left(a+\frac{-\alpha^2+\beta^2}{(7+4
n+\alpha+\beta) (8+4 n+\alpha+\beta)^2 (9+4
n+\alpha+\beta)}\right)$
\\[.8em]
\hline $\gamma_{n+1}^R$ & $\frac{32 (1+n) (3+2 n) (2+2 n+\alpha)
(3+2 n+\alpha) (2+2 n+\beta) (3+2 n+\beta) (2+2 n+\alpha+\beta)
(3+2 n+\alpha+\beta)} {(3+4 n+\alpha+\beta) (4+4 n+\alpha+\beta)^2
(5+4 n+\alpha+\beta)^2 (6+4 n+\alpha+\beta)^2 (7+4
n+\alpha+\beta)}$ \\[.5em]
\hline $\varrho_{n+1}^R$ & $p+\frac{(\alpha -\beta ) (\alpha+\beta )}{(4 n+\alpha +\beta +7)(4 n+\alpha +\beta +8)^2
(4n+\alpha +\beta +9)}+$\\
& $\frac{(\alpha -\beta ) (\alpha +\beta )}{(4n+\alpha +\beta+9) (4 n+\alpha +\beta +10)^2 (4 n+\alpha +\beta +11)}$
\\[.5em]
\hline $\rho_{n+1}^R$ &$4 (2 n+3) (2 n+\alpha +3) (2 n+\beta +3)(2 n+\alpha +\beta +3)\times \frac{\left(p+\frac{(\alpha -\beta)
(\alpha +\beta )}{(4 n+\alpha +\beta +5) (4 n+\alpha +\beta +6)^2(4 n+\alpha +\beta +7)}+ \frac{(\alpha -\beta ) (\alpha +\beta)}
{(4 n+\alpha +\beta +7) (4 n+\alpha +\beta +8)^2 (4 n+\alpha+\beta +9)}\right)}{(4 n+\alpha +\beta +5) (4n+\alpha +\beta +6)^2
(4 n+\alpha +\beta +7)}$
\\[.5em]
\hline \hline
all cases & $P_n$, $R_n$, $a_n$, $b_n$ nonorthogonal
\\
\hline
\end{tabular}
\end{center}

\begin{center}
\begin{tabular}[c]{p{2.5cm}||p{8.cm}}
 & {\bf GEGENBAUER}\\[.5em]
\hline \hline
$\beta_{n}\ \ | \ \ \gamma_{n+1}$
&  $0\ \ | \ \ \frac{(n+1)(n+2\alpha+1)}{(2n+2\alpha+1)(2n+2\alpha+3)}$    \\[.5em]
\hline \hline
$\beta_{n}^P$
& $\beta_0^P=ap+q+\frac 1{3+2\alpha}$ \\[.5em]
&$\beta_{n+1}^P=\frac{11+8 n^2 (1+2 q)+10 \alpha+4 n (1+2 q) (5+2
\alpha)+q (21+20 \alpha+4 \alpha^2)}
{21+16 n^2+20 \alpha+4 \alpha^2+8 n (5+2 \alpha)}$\\[.8em]
\hline $\gamma_{n+1}^P\ \ |\ \ \varrho_{n+1}^P$
&$\frac{4(1+n)(1+2 n)(1+n+\alpha)(1+2 n+2 \alpha)}{(3+4 n+2 \alpha)^2(5+16 n^2+12 \alpha+4 \alpha^2+8 n (3+2 \alpha))}\ \ |\ \ p$\\[.5em]
\hline $\rho_{n+1}^P$
& $\frac{4p (1+n)  (1+n+\alpha)}{15+16 n^2+16 \alpha+4 \alpha^2+16 n (2+\alpha)}$\\[.5em]
\hline \hline $\beta_{n}^R$
& $\beta_0^R=\frac{3-p^2 (5+2 \alpha)+q (5+2 \alpha)}{5+2\alpha}$\\[.5em]
& $\beta_{n+1}^R=\frac{23+8 n^2(1+2 q)+14\alpha+4 n(1+2 q) (7+2 \alpha)+q (45+28 \alpha+4 \alpha^2)}{45+16 n^2+28 \alpha+4 \alpha^2+8 n (7+2 \alpha)}$\\[.8em]
\hline $\gamma_{n+1}^R\ \ |\ \ \varrho_{n+1}^R$
& $\frac{4 (1+n) (3+2 n) (1+n+\alpha) (3+2 n+2 \alpha)}{(5+4 n+2 \alpha)^2 (21+16 n^2+20 \alpha+4 \alpha^2+8 n (5+2 \alpha))}\ \ |\ \ p$\\[.5em]
\hline $\rho_{n+1}^R$
& $p\frac{(3+2 n) (3+2 n+2 \alpha)}{35+16 n^2+24 \alpha+4 \alpha^2+16 n (3+\alpha)}$\\[.5em]
\hline \hline
$p=0$
&$P_n$, $R_n$ and $b_n=aR_n$ orthogonal  $a_n=0$  \\
\hline
$a=p=q=0$
&$P_n$ and $R_n$ orthogonal  $a_n=0$ and $b_n=0$ \\
\hline $p=q=0$
&$P_n$, $R_n$ and $b_n=aR_n$ orthogonal  $a_n=0$  \\
\hline $a=p=0$
&$P_n$ and $R_n$ orthogonal  $a_n=0$ and $b_n=0$ \\
\hline $a=q=0$
&$P_n$, $R_n$, $a_n$ and $b_n$ nonorthogonal   \\
\hline $a=0$
&$P_n$, $R_n$, $a_n$ and $b_n$ nonorthogonal \\
\hline $q=0$
&$P_n$, $R_n$, $a_n$ and $b_n$ nonorthogonal \\
\hline
\end{tabular}
\end{center}

\begin{center}
\begin{tabular}[c]{p{3.0cm}||p{7.5cm}}
  &  {\bf LEGENDRE} \\
\hline \hline
$\beta_{n}\ \ | \ \ \gamma_{n+1}$ & $0\ \ | \ \ \frac{4(n+1)^4}{(2n+1)(2n+2)^2(2n+3)}$     \\[.5em]
\hline \hline $\beta_{n}^P$ & $\beta_0^P=\frac 13+a p+q$\\
& $\beta_{n+1}^P=\frac{11+21 q+20 n (1+2 q)+8 n^2 (1+2 q)}{21+40n+16 n^2}$ \\[.5em]
\hline
$\gamma_{n+1}^P\ \ | \ \ \varrho_{n+1}^P\ \ |\ \ \rho_{n+1}^P$ & $\frac{4 (1+n)^2 (1+2 n)^2}{(3+4 n)^2(5+24 n+16 n^2)}\ \ | \ \ p\ \ |\ \ p\frac{4 (1+n)^2}{15+32 n+16 n^2}$
 \\[.5em]
\hline \hline
$\beta_{n}^R$ & $\beta_0^R=\frac 35-p^2+q$\\
& $\beta_{n+1}^R=\frac{23+45 q+28 n (1+2 q)+8 n^2 (1+2 q)}{45+56
n+16 n^2}$ \\[.5em]
\hline $\gamma_{n+1}^R\ \ | \ \  \varrho_{n+1}^R\ \ |\ \ \rho_{n+1}^R$ & $\frac{4 (1+n)^2 (3+2 n)^2}{(5+4 n)^2(21+40 n+16 n^2)} \ \ | \ \ p \ \ |\ \ p\frac{(3+2 n)^2}{35+48 n+16 n^2}$\\[.5em]
\hline \hline
$p=0$ & $P_n$, $R_n$, $b_n=aR_n$ orthogonal;  $a_n=0$ \\
\hline
$a=p=q=0$ & $P_n$, $R_n$ orthogonal;  $a_n=0$, $b_n=0$ \\
\hline
$p=q=0$ & $P_n$, $R_n$, $b_n=aR_n$ orthogonal;  $a_n=0$  \\
\hline
$a=p=0$ & $P_n$, $R_n$ orthogonal;  $a_n=0$, $b_n=0$ \\
\hline
$a=q=0$ &$P_n$, $R_n$, $a_n$, $b_n$ nonorthogonal  \\
\hline
$a=0$ & $P_n$, $R_n$, $a_n$, $b_n$ nonorthogonal \\
\hline
$q=0$ & $P_n$, $R_n$, $a_n$, $b_n$ nonorthogonal \\
\hline
\end{tabular}
\end{center}


\newpage

In the following tables, we consider the symmetric sequence  $\{Q_n\}_{n\geq 0}$, with recurrence coefficients
$\beta_n=0$ , $\gamma_{n+1}=\frac{1}{16}$ , $n\geq 0$.

\begin{center}
\begin{tabular}[c]{p{3.0cm}||p{7.5cm}}
 & \\
 & {\bf CHEBYSHEV 1$^{\mbox{\scriptsize{st}}}$ kind}  \\
\hline \hline
$\beta_{n} \ \ | \ \ \gamma_{n+1} $ & $0\ \ | \ \ \gamma_1=\frac 12\,,\:\:\:\gamma_{n+2}=\frac 14$     \\[.5em]
\hline \hline $\beta_{n}^P$ &
$\beta_0^P=\frac 12+a p+q\ ,\ \beta_{n+1}^P=\frac 12+q$ \\[.5em]
\hline
$\gamma_{n+1}^P\ \ | \ \ \varrho_{n+1}^P \ \ | \ \ \rho_{n+1}^P$ &
$\frac 1{16}\ \ | \ \ p \ \ | \ \ \frac p4 $
 \\[.5em]
\hline \hline $\beta_{n}^R$ &$\beta_0^R=\frac 34-p^2+q\ ,\ \beta_{n+1}^R=\frac 12+q$
\\[.5em]
\hline $\gamma_{n+1}^R\ \ | \ \  \varrho_{n+1}^R\ \ | \ \ \rho_{n+1}^R$ &
$\frac 1{16}\ \ | \ \ p  \ \ | \ \ \frac p4$ \\[.5em]
\hline \hline
$p=0\ ,\ q=-\frac{1}{2}$ & $R^{(k)}_n=P_n=Q_n$, $k\geq 1$\\
& $P_n(x)=A^{-n}W_{n}\Big(0;  \begin{array}{c}
0\\
\frac 12\\
\end{array}; 1 ; Ax
 \Big)\,,\:$ where $A=\pm 2$\\
\hline \hline
& \\
& {\bf CHEBYSHEV 2$^{\mbox{\scriptsize{nd}}}$ kind } \\
\hline \hline
$\beta_{n}\ \ | \ \ \gamma_{n+1}$ &  $0\ \ | \ \ \frac{1}{4}$    \\[.5em]
\hline \hline $\beta_{n}^P$ &$\beta_0^P=\frac 14+a p+q\ ,\ \beta_{n+1}^P=\frac 12+q$     \\[.5em]
\hline $\gamma_{n+1}^P\ \ | \ \ \varrho_{n+1}^P \ \ | \ \ \rho_{n+1}^P$ &
$\frac 1{16} \ \ | \ \  p  \ \ | \ \  \frac p4$
 \\[.5em]
\hline \hline $\beta_{n}^R$ &$\beta_0^R=\frac 12-p^2+q\ ,\ \beta_{n+1}^R=\frac 12+q$ \\[.5em]
\hline $\gamma_{n+1}^R\ \ | \ \  \varrho_{n+1}^R\ \ | \ \ \rho_{n+1}^R$ &
$\frac 1{16} \ \ | \ \ p  \ \ | \ \ \frac p4$  \\[.5em]
\hline \hline
$p=0\ ,\ q=-\frac{1}{2}$ &  $P^{(k)}_n=R_n=Q_n$, $k\geq 1$ ; $R_n(x)=A^{-n}W_n(Ax)\,,\:$ where $A=\pm 2$\\
\hline \hline
&\\
& {\bf CHEBYSHEV 1$^{\mbox{\scriptsize{st}}}$ and 2$^{\mbox{\scriptsize{nd}}}$ kinds} \\ \hline \hline
$p=0$ & $P_n$, $R_n$, $b_n=aR_n$ orthogonal;  $a_n=0$ \\
& $P_n$, $R_n$ positive definite, if $a$, $q \in {\bkR}$\\
& $P^{(k)}_n=R^{(k)}_n\ ,\ n\geq 1\ ,\ k\geq 1$\\
& $P_n(x)=R_n(-\frac{1}{4};x)$ \\
\hline
$p=0\ ,\ q=-\frac{1}{2}$ & $P_n(x)=R_n(-\frac{1}{4};x)$
\\
\hline \hline
$a=p=q=0$ & $P_n$, $R_n$ orthogonal;  $a_n=0$, $b_n=0$  \\
\hline
$p=q=0$ &$P_n$, $R_n$, $b_n=aR_n$ orthogonal;  $a_n=0$  \\
\hline
$a=p=0$ & $P_n$, $R_n$ orthogonal;  $a_n=0$, $b_n=0$ \\
\hline
$a=q=0$ & $P_n$, $R_n$, $a_n$, $b_n$ nonorthogonal  \\
\hline
$a=0$ & $P_n$, $R_n$, $a_n$, $b_n$ nonorthogonal \\
\hline
$q=0$ & $P_n$, $R_n$, $a_n$, $b_n$ nonorthogonal\\
\hline
\end{tabular}
\end{center}

\begin{center}
\begin{tabular}[c]{p{3.0cm}||p{7.5cm}}
& \\
 & {\bf CHEBYSHEV  3$^{\mbox{\scriptsize{rd}}}$ kind }\\
\hline \hline
$\beta_{n}\ \ |\ \ \gamma_{n+1} $ & $\beta_0=\frac 12\ ,\ \beta_{n+1}=0 \ \ |\ \ \frac 14$ \\[.5em]
\hline \hline $\beta_{n}^P$ & $\beta_0^P=a
\left(p+\frac{1}{2}\right)+q+\frac{1}{4}\ \ |\ \ \beta_{n+1}^P=\frac 12+q$ \\[.5em]
\hline \hline $\gamma_{n+1}^P\ \ | \ \ \varrho_{n+1}^P \ \ | \ \ \rho_{n+1}^P$ &
$\frac 1{16}\ \ |\ \  p  \ \ |\ \  \frac p4$ \\[.5em]
\hline \hline $\beta_{n}^R$ & $\beta_0^R=\frac 12-p\left(\frac
12+p\right)+q\ ,\   \beta_{n+1}^R=\frac 12+q$\\[.5em]
\hline $\gamma_{n+1}^R\ \ | \ \  \varrho_{n+1}^R\ \ | \ \ \rho_{n+1}^R$ &
$\frac 1{16}\ \ | \ \ p \ \ | \ \ \frac p4$ \\[.5em]

\hline \hline
$p=0$ & $P_n(x)=R_n(\frac{1}{2}(a-\frac 12);x)$ \\
\hline
$p=0\ ,\ q=-\frac{1}{2}$ & $P_n(x)=R_n(\frac{1}{2}(a-\frac 12);x)$  \\
& $R_n(x)=A^{-n}W_n\left(-\frac 12;Ax\right)\,,\:$ where $A=\pm 2$\\
\hline \hline

& \\
& {\bf CHEBYSHEV 4$^{\mbox{\scriptsize{th}}}$ kind}  \\
\hline \hline
 $\beta_{n}\ \ |\ \ \gamma_{n+1}$ &$\beta_0=-\frac 12\ ,\ \beta_{n+1}=0\ \ |\ \ \frac{1}{4}$
 \\[.5em]

\hline \hline $\beta_{n}^P$
&$\beta_0^P=a\left(p-\frac{1}{2}\right)+q+\frac{1}{4}\ ,\ \beta_{n+1}^P=\frac 12+q$\\[.5em]

\hline $\gamma_{n+1}^P\ \ | \ \ \varrho_{n+1}^P \ \ | \ \ \rho_{n+1}^P$ &
$\frac 1{16}\ \ |\ \  p  \ \ |\ \  \frac p4$
\\[.5em]

\hline \hline $\beta_{n}^R$  & $\beta_0^R=\frac{1}{2} \left(-2 p^2+p+2
q+1\right)\ ,\ \beta_{n+1}^R=\frac 12+q$
\\[.5em]

\hline $\gamma_{n+1}^R\ \ | \ \ \varrho_{n+1}^R \ \ | \ \ \rho_{n+1}^R$ &
$\frac 1{16}\ \ |\ \  p  \ \ |\ \  \frac p4$
\\[.5em]

\hline \hline
$p=0$ &  $P_n(x)=R_n(-\frac{1}{2}(a+\frac{1}{2});x)$\\
\hline
$p=0\ ,\ q=-\frac{1}{2}$ &  $P_n(x)=R_n(-\frac{1}{2}(a+\frac{1}{2});x)$\\
& $R_n(x)=A^{-n}W_n\left(\frac 12;Ax\right)\,,\:$ where $A=\pm 2$\\
\hline \hline
& \\
& {\bf CHEBYSHEV 3$^{\mbox{\scriptsize{rd}}}$ and 4$^{\mbox{\scriptsize{th}}}$ kinds } \\
\hline \hline
$p=0$ & $P_n$, $R_n$, $a_n$, $b_n$ orthogonal\\
& $P_n$, $R_n$ positive definite, if $a$, $q \in {\bkR}$\\
& $P^{(k)}_n=R^{(k)}_n\ ,\ n\geq 1\ ,\ k\geq 1$\\
 \hline
 $p=0\ ,\ q=-\frac{1}{2}$ & $P_n$, $R_n$, $a_n$, $b_n$ orthogonal\\
&  $P^{(k)}_n=R_n=Q_n$, $k\geq 1$\\
 \hline \hline
$a=p=q=0$ & $P_n$,
$R_n$, $a_n$, $b_n$ orthogonal\\
\hline $p=q=0$ &$P_n$,
$R_n$, $a_n$, $b_n$ orthogonal
\\
\hline $a=p=0$ & $P_n$,
$R_n$, $a_n$, $b_n$ orthogonal
\\
\hline $a=q=0$ &
$P_n$, $R_n$, $a_n$, $b_n$ nonorthogonal
\\
\hline $a=0$ &
$P_n$, $R_n$, $a_n$, $b_n$ nonorthogonal
\\
\hline $q=0$ &
$P_n$, $R_n$, $a_n$, $b_n$ nonorthogonal
\\
\hline
\end{tabular}
\end{center}




\end{document}